\documentclass[arxiv,reqno,twoside,a4paper,12pt]{amsart}
\usepackage{mltools}

\withdetails
\usepackage{mlmath62,mlthm10-REAR}

\usepackage{upref,amsmath,amssymb,amsthm,mathrsfs}
\usepackage{pgf,tikz, graphics, xypic, latexsym}
\usepackage{graphicx}
\usepackage[colorlinks=true,linkcolor=blue,citecolor=blue]{hyperref}

\usepackage[scaled]{helvet} 
\usepackage{courier} 
\usepackage[mathbf]{euler}
\usepackage{mllocal}

\numberwithin{equation}{section}

\definecolor{qqwuqq}{rgb}{0,0,0}

\begin{document}

\title[Newton Polyhedra and Whitney Equisingularity for IDS]
{Newton Polyhedra and Whitney Equisingularity for Isolated Determinantal Singularities}

\author{Tha\'is M. Dalbelo}
\address{Universidade Federal de S\~{a}o Carlos, Brazil}
\email{thaisdalbelo@ufscar.br}

\author{Luiz Hartmann}
\address{Universidade Federal de S\~{a}o Carlos, Brazil}
\email{hartmann@dm.ufscar.br, luizhartmann@ufscar.br}
\urladdr{http://www.dm.ufscar.br/profs/hartmann}

\author{Maicom Varella}
\address{Universidade Federal de S\~{a}o Carlos, Brazil}
\email{maicomvarella@estudante.ufscar.br}

\thanks{ During this work T. M. Dalbelo was partially supported by FAPESP-Grant
2019/21181-0 and by FAPESP-Grant 2023/01018-2. L. Hartmann was partially supported by FAPESP-Grant
2022/16455-6.  M. Varella was supported by CAPES-Brazil Finance code 001 and by
DAAD-Germany Doctoral Research Grant (Number 57507871)}

\date{\today}

\begin{abstract}

Using Newton polyhedra and non-degeneracy of matrices we present conditions
which guarantee the Whitney equisingularity of families of isolated
determinantal singularities.

\end{abstract}



\maketitle

\tableofcontents


\section{Introduction}

Newton polyhedra are an important tool in the study of topological, geometric and algebraic properties of analytical varieties and their deformations. In Singularity Theory there are countless articles, where Newton polyhedra and non-degeneracy conditions are applied in the study of hypersurfaces, complete intersection and toric varieties, for instance.

For determinantal varieties, the notion of Newton polyhedra and the non-degeneracy conditions were studied by Esterov \cite{Esterov}. In this article, the author introduced a concept of non-degenerate matrices and applied it to compute the multiplicity of a non-degenerate determinantal variety and the Euler characteristic of the Milnor fiber of a function defined on it, both in terms of Newton polyhedra.

	 The concept of Newton polyhedra is strongly related to the Whitney equisingularity of families of varieties. On unpublished notes Brian\c con studied the Whitney equisingularity for a family of Newton non-degenerate isolated hypersurface singularities. Eyral and Oka \cite{EyralOka} presented a criterion to verify the	Whitney equisingularity for families of (possibly non-isolated) non-degenerate hypersurface singularities.

	 	  On the other hand the concept of Whitney equisingularity is also strongly related to the polar multiplicities for many important classes of spaces. For instance, in \cite{Gaffney, Gaffney2}, Gaffney showed that a family of d-dimensional isolated complete intersection singularities (ICIS) $\left\{  (X_t,0) \right\}_{t\in D}$ of any dimension $d$ is Whitney equisingular if, and only if, the polar multiplicities $m_i(X_t, 0)$, $i = 0,\dots, d$ are constant on this family, where $D$ is an open disc around the origin in $\mathbb{C}$.

	 	The polar multiplicities $m_i(X,0)$, for $i = 0,\dots, d-1$ are defined for any variety $(X, 0)$ (see \cite{LT}), while $m_d(X, 0)$ was defined, initially, only for ICIS in \cite{Gaffney}. Using topological and geometric information from generic linear projections Nu\~no-Ballesteros, Or\'efice-Okamoto and Tomazella \cite{BBT} defined the top polar multiplicity $m_d(X_A^s, 0)$ for an isolated determinantal singularity (IDS) $(X_A^s,0)$.  Moreover, the authors proved (see \cite[Theorem 4.3]{BBT}), that the multiplicity $m_d(X_A^s,0)$ is related to the vanishing Euler characteristic of $X_A^s$ (see Definition \ref{vanishing} below) by the following Lê-Greuel type formula
	 	\begin{equation}\label{toppolar2}
	 	m_d (X_A^s,0) = \nu(X_A^s, 0) + \nu(X_A^s\cap p^{-1}(0),0),
	 	\end{equation}
	 	where $p: (\mathbb{C}^m,0) \to (\mathbb{C},0)$ is a generic linear function and $d = \dim X_A^s$. For codimension $2$ IDS in $\mathbb{C}^4$ and $\mathbb{C}^5$ the multiplicity $m_d(X_{A}^s,0)$ was also defined by Pereira and Ruas in \cite{MC} and the authors also proved the Lê-Greuel type formula \Eqref{toppolar2} (see \cite[Theorem 5.2]{MC}).

	 	In \cite[Theorem 5.3]{BBT4}, the authors proved that a family of IDS $\left\{(X_{A_t}^s,0) \right\}$ is Whitney equisingular if, and only if, $m_i(X_{A_t}^s, 0)$, $i = 0,\dots, d$, do not depend on $t \in D$.

	The main goal of this work is to use Newton polyhedra and non-de\-ge\-ne\-ra\-cy of matrices to present a criterion to test the Whitney equisingularity for families of IDS, which is the following theorem.

	\begin{theorem}
		\label{equisingularity}
		Let $\left\{(X_{A_t}^s, 0) \right\}_{ t\in D }$, be a d-dimensional family of determinantal singularities, defined by the germ of matrices $A_t = ((a_{ij})_t):(\mathbb{C}^m,0)\rightarrow (M_{n,k},0)$ with holomorphic entries. Suppose that $X_{A_0}^s$ has an isolated singularity at $0$ and, for all $t \in D$, the matrix $A_t$ satisfies the following conditions:
		\begin{enumerate}
			\item the matrix $A_t$ is Newton non-degenerate (Definition \ref{non-degeneracity} below);
			\item the Newton polyhedra $\Delta^t_{j}$ of $(a_{ij})_t$ are convenient and independent of $t$.

		\end{enumerate}
		Then the family $\left\{(X_{A_t}^s, 0) \right\}_{ t\in D}$ is Whitney equisingular.
	\end{theorem}

	In order to prove Theorem \ref{equisingularity} we developed the following steps. Firstly, we apply the notions introduced by Esterov \cite{Esterov} and the ideas used by Eyral and Oka \cite{EyralOka} to prove (under the assumptions of Theorem \ref{equisingularity}) the existence of a positive number $\varepsilon > 0$ for which the set $(X_{A_t}^s\setminus\{0\}) \cap B_{\varepsilon}$ is non-singular and $B_{\varepsilon}$ is a Milnor ball for $(X_{A_t}^s,0)$ independent of $t \in D$, this is proved in Corollary \ref{cororadius}. This result is in fact more general and holds without the assumption of the Newton polyhedra being convenient, which is proved in Proposition \ref{radius}. More precisely Proposition \ref{radius} holds for determinantal varieties with arbitrarily singular set, not necessarily the ones with zero dimensional singular set.  As an important consequence of these results we obtain in Corollary \ref{vanishing-constant}, the constancy of vanishing Euler characteristic for the family $\{(X_{A_t}^s,0)\}$.

	 In Proposition \ref{radius-section} and Corollary \ref{cororadius-section} we prove a version of Proposition \ref{radius} and Corollary \ref{cororadius}, respectively, to the family $\{(X_{A_t}^s \cap p^{-1}(0),0)\}$. Finally, we use these results to prove Theorem \ref{equisingularity}.

		We remark that Whitney equisingularity is a very strong form of equisingularity. A Whitney equinsigular family of varieties $\{(X_t,0)\}_{t \in D}$ is also topologically equisingular, by the Thom-Mather first isotopy theorem (\cite{Mather},\cite{Thom}). This means that the local, ambient, topological type of $X_t$ at $0$ does not depend of $t$ for $t\in D$.

The results mentioned above, are state and proved in Section $3$. In the following section we present some definitions and properties regarding isolated determinantal singularities and Newton polyhedra.

\section{Determinantal deformation and Newton polyhedra}

\subsection{Determinantal deformation}

Determinantal varieties have been wi\-dely studied by researchers in
Commutative Algebra, Algebraic Geometry and   Singularity Theory. Just to mention a few  works, we can quote, for instance, Bruns and Herzog \cite{Bruns2}, Bruns and Vetter \cite{Bruns1}, Ebeling and Gusein-Zade \cite{EGZ}, Frühbis-Krüger
and Neumer \cite{AFG2}, Gaffney, Grulha Jr. and Ruas
\cite{GGR}, Nuño-Ballesteros, Oréfice-Okamoto and Tomazella
\cite{BBT} and Pereira and Ruas \cite{MC}.

Consider $M_{n,k}$ the set of complex matrices of size $n\times k$ and
$M_{n,k}^s$ the subset consisting of the matrices with rank less than $s$,
where $0< s\leq n \leq k$ are integers. The set $M_{n,k}^s$ is an irreducible
subvariety of $M_{n,k}$ with codimension $(n-s+1)(k-s+1)$, which is called
{\it generic determinantal variety}, and $M_{n,k}^{s-1}$ is its singular set.

Let $A:(\mathbb{C}^m,0) \rightarrow (M_{n,k},0)$ be a holomorphic map germ defined by $A(x)=(a_{ij}(x))$ with $a_{ij} \in \mathcal{O}_m$, where $\mathcal{O}_m$ is the ring of holomorphic functions in $\mathbb{C}^m$, for $1 \leq i \leq n$ and $1 \leq j \leq k$.

\begin{definition}
Let $(X_A^s,0) \subset (\mathbb{C}^m,0)$ be the germ of the variety defined by $X_A^s=A^{-1}(M_{n,k}^s)$. We say that $X_A^s$ is a germ of {\it determinantal variety} of type $(n,k;s)$ in $(\mathbb{C}^m,0)$, if its dimension is equal to $$ m- (n-s+1)(k-s+1).$$
\end{definition}

The analytical structure of $X_A^s$ is the one defined by $A$ and $M_{n,k}^s$, \ie it is given by the $s$ size minors of $A(x)$ (not necessarily reduced). When $s=1$, the determinantal variety $X_A^s$ is a complete intersection variety.
We observe that the singular locus of $X_A^s$ is defined via the Jacobian matrix of
the defining equations. In particular, the singular locus contains all the non
reduced components of the determinantal germ (if any). For more details see \cite[pg 31]{Hartshorne}.

In this work, we will use the concept of isolated determinantal singularity (IDS), which was defined by
 Nu\~no-Ballesteros, Oréfice-Okamoto and Tomazella in \cite{BBT}.




\begin{definition}
Let $(X_A^s,0) \subset (\mathbb{C}^m,0)$ be a determinantal variety of type $(n,k;s)$.
The variety $(X_A^s,0)$ is said to be an {\it isolated determinantal singularity} (IDS) if
it satisfies the condition $$s=1 \textrm{ or } m<(n-s+2)(k-s+2),$$
 $X_A^s$ is smooth at $x\in X^s_A$ and $\textrm{rank } A(x)=s-1$ for all $0\neq x\in X^s_A$ in a neighborhood of the origin.
\end{definition}


 We observe that if $(X_A^s, 0)$ is an IDS of dimension greater than $0$, then it is reduced, because it is
Cohen–Macaulay and satisfies Serre’s $R0$ condition (see, for instance, \cite{GP}).

\begin{remark}
The theory of singularities for matrices has been established by several authors. For instance, Bruce and Tari utilized $\mathcal{G}$-equivalence of matrices to classify simple germs of families of symmetric, anti-symmetric, and quadratic matrices \cite{Bruce, BruceTari}, while Frühbis-Krüger and Neumer \cite{AFG1, AFG2} focused on the case of $n \times (n+1)$ matrices. Pereira \cite{MP} studied the $n \times k$ case. The author proved that isolated determinantal singularities are $\mathcal{G}$-determined, implying they have a polynomial representative. Here, $\mathcal{G} = \mathcal{R} \ltimes \operatorname{GL}_k(\mathcal{O}_m) \times \operatorname{GL}_n(\mathcal{O}_m)$, where $\mathcal{R}$ represents the coordinate change group in $(\mathbb{C}^m,0)$, and $\operatorname{GL}_i(\mathcal{O}_m)$ denotes the group of invertible $i \times i$ matrices with entries in $\mathcal{O}_m$.
\end{remark}

In the following, we present the definition of vanishing Euler characteristic.

\begin{definition}\label{vanishing}
	The {\it vanishing Euler characteristic} of an IDS $X_A^s$, denoted by
	$\nu(X_A^s,0)$, is defined as $$\nu(X_A^s, 0) = (-1)^d(\chi(\tilde{X_{A}^s}) - 1),$$
	where $\tilde{X_{A}^s}$ is the generic fiber of the determinantal smoothing\footnote{For the definition see \cite[Definition 3.3]{BBT}.} of $X_A^s$, $\chi( \cdot )$ denotes the Euler characteristic and $d=\dim X_A^s$.
\end{definition}

For more details on the vanishing Euler characteristic and the Milnor number of IDS see \cite[Definition 3.2]{BBT} or \cite[Definition 5.1]{MC}.

Now, we recall  some notions about determinantal deformations. Consider a map germ
$\mathcal{A}: (\mathbb{C}^m \times \mathbb{C}, 0) \to (M_{n,k},0)$ such that
$\mathcal{A}(x, 0) = A(x)$ for all $x \in \mathbb{C}^m$. When $X_{A}^s$ is a determinantal variety, the projection
$$\begin{array}{cccc}
\pi: &(X_{\mathcal{A}}^s,0) & \to & (\mathbb{C},0) \\
     & (x,t) & \mapsto & t
\end{array}$$
is called a {\it determinantal deformation of} $X_{A}^s$ . If we fix a small enough representative $A :
B_{\varepsilon} \to (M_{n,k},0)$, where $B_{\varepsilon}$ is the open ball
centred
	at the origin with radius $\varepsilon > 0$, then we set $A_t(x) := \mathcal{A}(x, t)$ and $X_{A_t}^s = A^{-1}_t (M^{s}_{n,k})$.

If $X_{A_t}^s$ is a determinantal deformation of $(X_{A_0}^s, 0)$ as above, we say that:
\begin{enumerate}
	\item $X_{A_t}^s$ is origin preserving if $0 \in S(X_{A_t}^s)$, for all $t$ in $D$, where $S(X_{A_t}^s)$ denotes the singular
	set of $X_{A_t}^s$ and $D \subset \mathbb{C}$ is a disc around the origin. Then $\left\{(X_{A_t}^s, 0) \right\}_{ t\in D }$ is called a {\it $1$-parameter family} of IDS;
	\item $\left\{(X_{A_t}^s, 0) \right\}_{ t\in D }$ is a good family if there exists $\varepsilon > 0$ with $S(X_{A_t}^s) = \left\{0\right\}$ on $B_{\varepsilon}$, for all $t$ in $D$;
	\item $\left\{(X_{A_t}^s, 0) \right\}_{ t\in D }$ is Whitney equisingular if it is a good family and $\left\{X_{\mathcal{A}}^s \setminus T, T \right\}$ satisfies the Whitney
	conditions, where $T = \left\{0\right\} \times D$.
\end{enumerate}

In \cite[Theorem 5.3]{BBT4}, the authors proved that a family of IDS $\left\{(X_{A_t}^s,0) \right\}_{t\in D}$ is Whitney equisingular if, and only if, $m_i(X_{A_t}^s, 0)$, $i = 0,\dots, d$, do not depend on $t \in D$, where $D$ is an open disc around the origin in $\mathbb{C}$. Moreover, they also obtained the following result about the constancy of the vanishing Euler characteristic \cite[Corollary 4.3]{BBT4}.

\begin{corollary}\label{corovanishing}
	Let $X_{A_t}^s$ be a determinantal deformation of the IDS $X_{A_0}^s$. Then, the sum $\sum_{x\in S(X_{A_t}^s)}\nu (X_{A_t}^s,x)$ is constant on $t$ if and only if $\chi (X_{A_t}^s) =1$, where $S(\cdot)$ denotes the singular set.
\end{corollary}

\subsection{Newton polyhedra}

 The Newton polyhedra of polynomial functions are important objects which can be very useful to compute some invariants, such as the Milnor number, local Euler obstruction, multiplicities, among others.

The monomial $x_1^{a_1} \cdots  x_m^{a_m}$ is denoted by $x^a$, where $a=(a_1, \dots , a_m) \in \mathbb{Z}^m$. We denote by $\mathbb{R}^m_{+}$, the positive orthant of $\mathbb{R}^m$. A subset $\Delta \subset \mathbb{R}^m_{+}$ is called a {\it Newton
polyhedron} when there exists some $P\subset \mathbb{Z}^m_{+}$ such that
$\Delta$ is the convex hull of the set $\{p+v: p\in P \textrm{ and } v \in
\mathbb{R}^m_{+}\}$. In this case, $\Delta$ is said to be the Newton
polyhedron determined by $P$.

\begin{definition}
If the polyhedron $\Delta$ touches all the coordinate axes, we say that $\Delta$ is {\it convenient}.
\end{definition}

\begin{definition}
If $f\in \mathcal{O}_m$ is a germ of a
polynomial function $f(x)=\displaystyle \sum_{a \in \mathbb{Z}^m_{+}} c_a x^a$,
then the {\it support} of $f$ is $\rm{supp}(f):=\setdef{a\in
\mathbb{Z}^m_{+}}{ c_a \neq 0}$. The {\it Newton polyhedron} of $f$ $\Delta_f$ is
the Newton polyhedron determined by $\rm{supp}(f)$.
\end{definition}

Let $P=(p_1,\dots ,p_m)$ be a weight vector. For each $x\in \mathbb{R}^m$, we define $P(x)=\sum_{i=1}^mp_i\cdot x_i$. For a positive weight ($p_i \geq 0$, for $i=1, \dots, m$), we define $d(P;f)$ as the minimal value of the restriction $P|_{\Delta_f}$ \ie $d(P;f)=\min\{P(x): x\in \Delta_f\}$. Let $$\Gamma(P;f) = \setdef{ x\in \Delta_f}{ P(x) = d(P;f)}$$ be the unique face of $\Delta_f$, where $P$ takes the minimal value $d(P;f)$. For a strictly positive weight $P$ ($p_i > 0$, for $i=1, \dots, m$), $\Gamma(P;f)$ is a bounded face of $\Delta_f$ (it can even be a $0$-dimensional face, i.e. a vertex of $\Delta_f$).  When there is no risk of misunderstandings, we denote $\Gamma(P;f)$ simply by $\Gamma$. We define $$f|_{\Gamma}(x) = \sum_{a\in \Gamma} c_ax^a$$ and we call $f|_{\Gamma}$ the face function of $f$ with respect to $\Gamma$.

In the following, we present the Newton non-degeneracy conditions adapted to
determinantal singularities. The Newton non-degeneracy conditions in this case were first defined in \cite[Definition 1.16]{Esterov} for determinantal singularities. The definition presented by Esterov comprehends determinantal singularities given by the maximal minors of a matrix $A(x) = (a_{ij}(x))$, under the assumption that, all the entries of a column $j$ has the same Newton polyhedron $\Delta_j$. Here we still consider that, for each $j=1,\dots ,k$, the Newton polyhedron of $a_{i,j}$ is $\Delta_j$, for all $i=1,\dots ,n$. However, we do not assume that the variety is given only by the maximal minors.

In the next definition we use the Minkowski sum between sets. More precisely, given $A,B \subset \mathbb{R}^m$ the Minkowski sum is defined by
 $$
 A+B = \setdef{x+y}{  x \in A, y \in B }.
 $$
We also remark that, in the next definition, the Newton polyhedra are not necessarily convenient and the  considered faces can be of arbitrarily dimension.

\begin{definition}\label{non-degeneracity}
Let $A = (a_{ij}): (\mathbb{C}^m,0) \to (M_{n,k},0)$ be a germ of a holomorphic matrix. We denote by $\Delta_j \subset \mathbb{R}^m_{+}$ the Newton polyhedron of $a_{ij}$, for all $i=1,\dots,n$.
\begin{enumerate}
\item The matrix $A$ is said to be {\it Newton non-degenerate} if, for each collection of faces  $\Gamma_j \subset \Delta_j$ such that the sum $\sum_{j=1}^k \Gamma_j$ is a bounded face of the polyhedron $\sum_{j=1}^k \Delta_j$, the
	polynomial matrix $(a_{ij}|_{\Gamma_j})$ defines a non-singular
	determinantal variety of type $(n,k;s)$ in $(\mathbb{C}\setminus \{0\})^m$ or the empty set inside $(\mathbb{C}\setminus \{0\})^m$. In this case, we say that $X_A^s$ is a Newton non-degenerate determinantal singularity.

\label{non-degeneracity-variety}
\item Let $f_{1},\dots,f_l :(\mathbb{C}^m,0)\to (\mathbb{C},0)$ be holomorphic functions with Newton polyhedra $\Delta'_{1},\dots ,\Delta_l'$, respectively. The variety $X_A^s\cap f_{1}^{-1}(0)\cap \cdots \cap f_l^{-1}(0)$ is said to be {\it Newton non-degenerate} if, $A$ is a Newton non-degenerate matrix and for each collection of faces  $\Gamma_j \subset \Delta_j$, $j=1,\dots ,k$ and $\Gamma'_{\gamma}\subset \Delta'_{\gamma}$,  $\gamma = 1,\dots ,l$ such that the sum $\sum_{j=1}^k \Gamma_j+\sum_{\gamma =1}^l \Gamma'_{\gamma}$ is a bounded face of the polyhedron $\sum_{j=1}^k \Delta_j+\sum_{\gamma =1}^l \Delta'_{\gamma}$, the variety $$X_{(a_{ij}|_{\Gamma_j})}^s\cap (f_{1}|_{\Gamma'_{1}})^{-1}(0)\cap \cdots \cap (f_l|_{\Gamma'_{l}})^{-1}(0)$$
	is a non-singular
	variety in $(\mathbb{C}\setminus \{0\})^m$ with dimension equal to  $\dim X_A^s -l$ or the empty set inside $(\mathbb{C}\setminus \{0\})^m$.
\end{enumerate}
\end{definition}

When $l=1$, the definition introduced by Esterov (\cite[Definition 1.16]{Esterov}) of $f$ being non-degenerate with respect to a matrix $A$ implies that the variety $X_A^s\cap f^{-1}(0)$ is a Newton non-degenerate variety.

\begin{remark}
Given a Newton non-degenerate matrix $A$ and a generic linear form  with respect to $X_A^s$ $h: (\mathbb{C}^m,0) \to (\mathbb{C},0)$, the restriction of $h$ to $X_A^s$ may be degenerate if we eliminate one variable using $h=0$ (see \cite[Example (I-2)]{Oka}).
However the variety $X_A^s\cap h^{-1}(0) \subset \mathbb{C}^m$ is Newton non-degenerate.
\end{remark}

\section{Whitney equisingularity}

The main purpose of this section is to prove Theorem \ref{equisingularity}. Before proving it, we need some concepts and results in order to guarantee an uniformity condition on the family $\{(X_{A_t}^s,0)\}$. The following concept was first introduced by Eyral and Oka \cite{EyralOka,Oka} for hypersurfaces and complete intersection singularities and we extend it to determinantal singularities.

\begin{definition}
Let $f:(\mathbb{C}^m,0)\to (\mathbb{C},0)$ be a germ of an analytic function and $I\subset \{1,\dots ,m\}$. We say that $\mathbb{C}^I=\{(x_1,\dots,x_m)\in \mathbb{C}^m: x_i = 0\; ,i \notin I\}$ is an {\it admissible coordinate subspace for} $f$ if $f^I:=f|_{\mathbb{C}^I\cap U}$ is not constantly zero, where $U$ is a neighbourhood of the origin in $\mathbb{C}^m$.
\end{definition}

\begin{definition}
Let $A=(a_{ij}):(\mathbb{C}^m,0) \to (M_{n,k},0)$ be a germ of a holomorphic map, we say that $\mathbb{C}^I$ is an {\it admissible coordinate space for} $A$ if $\mathbb{C}^I$ is admissible for each $a_{ij}$, $i = 1, \dots ,n$ and $j=1,\dots ,k$ .
\end{definition}

Along this work, we will use the following notation:
\begin{enumerate}
\item $A^I = (a_{ij}^I)$,  where $a_{ij}^{I} = a_{ij}|_{\mathbb{C}^{I}\cap U}$;
\item $A^{*I} = (a_{ij}^{*I})$, where $a_{ij}^{*I} = a_{ij}|_{\mathbb{C}^{*I}\cap U}$ and $\mathbb{C}^{*I}=\{(x_1,\dots ,x_m)\in \mathbb{C}^m: x_i=0 \Leftrightarrow i \notin I\}$.
\end{enumerate}

In the following Lemma we extend \cite[Proposition 7]{Oka} to determinantal singularities and using the ideas of Eyral and Oka \cite[Proposition 3.1]{EyralOka}, we prove Proposition \ref{radius}.

 \begin{lemma}
\label{lemmaadmissible}
Let $A:(\mathbb{C}^m,0)\to (M_{n,k},0)$ be a germ of a holomorphic matrix and $f_{1},\dots ,f_l: (\mathbb{C}^m,0)\to (\mathbb{C},0)$ be germs of functions. Let $\mathbb{C}^I$ be an admissible coordinate subspace for $A,f_{1},\dots f_l$.
\begin{enumerate}
\item If the matrix $A$ is Newton non-degenerate, then $A^I$ is Newton non-degenerate.
\item If the variety $X_A^s\cap f_{1}^{-1}(0)\cap \cdots \cap f_l^{-1}(0)$ is Newton non-degenerate, then $X_{A^I}^s\cap (f_{1}^I)^{-1}(0)\cap \cdots \cap (f_l^I)^{-1}(0)$ is a Newton non-degenerate variety.
\end{enumerate}
\end{lemma}

\begin{proof}
Since booth cases are analagous, we prove only $(i)$.
For each $j=1,\dots,k$ let $\Gamma_j$ be the faces of $\Delta_{a_{ij}^I}$ such that $\sum_{j=1}^k \Gamma_j$ is a bounded face of $\sum_{j=1}^k \Delta_{a_{ij}^I}$. Since $\Delta_{a_{ij}^I} = \Delta_{a_{ij}}\cap \mathbb{R}^I$, $\Gamma_j$ is also a face of $\Delta_{a_{ij}}$ such that $\sum_{j=1}^k \Gamma_j$ is a bounded face of $\sum \Delta_{a_{ij}}$ for each $j=1,\dots,k$. Since the matrix $A$ is Newton non-degenerate, then $(a_{ij}^I|_{\Gamma_j})$ defines a non-singular determinantal variety in $(\mathbb{C}\setminus \{0\})^m$. Hence, $A^I$ is Newton non-degenerate.
\end{proof}

In \cite[Proposition 3.1]{EyralOka} Eyral and Oka proved an uniform version of \cite[ChapterIII, Lemma(2.2)]{Oka} and \cite[Theorem 19]{Oka3}. Here we extend \cite[Proposition 3.1]{EyralOka} to families of non-degenerate determinantal singularities.

 	\begin{proposition}
	\label{radius}
	Suppose that for all t sufficiently small, the following two conditions are satisfied:
	\begin{enumerate}
		\item the matrix $A_t = ((a_{ij})_t)$ is Newton non-degenerate;
		\item the Newton polyhedron $\Delta_j^t$ of $(a_{ij})_t$ is independent of $t$, for all $j=1,\dots ,k$.
	\end{enumerate}
	Then there exists a positive number $R > 0$ such that for any admissible coordinate subspace $\mathbb{C}^I$ of $A_0$ and any $t$ sufficiently small, the set $X_{A_t}^s\cap \mathbb{C}^{*I} \cap B_R$ is non-singular and intersects transversely with $S_r$ for any $r<R$, where $B_R$ (respectively, $S_r$) is the open ball (respectively, the sphere) with center at the origin $0 \in \mathbb{C}^m$ and radius $R$ (respectively, $r$).
\end{proposition}

\begin{proof} This proof is divided in 4 steps. Firstly, we prove the smoothness part, which is splitted in two steps, Claim $1$ and Claim $2$ below. The transversality part will be proved using Claim $3$ and Claim $4$.

 In the sequence we introduce the objects which we need to state Claim $1$ and $2$.

As $\Delta_j^t$ is independent of $t$, the set of admissible systems is also independent of $t$ and, since there are only finitely many subsets $I \subset \{1,\dots,m\}$, it suffices to show the result for a fixed $I=\{1,\dots ,r\}$, $r\leq m$.

 We argue by contradiction. Suppose that for all $R>0$ the intersection $X_{A_t}^s\cap \mathbb{C}^{*I}\cap B_R$ has a singular point. Consider the sequence $\{(t_{R},z_{R})\}$ of points in $X_{\mathcal{A}}^s\cap(D\times \mathbb{C}^{*I})$ converging to $(0,0)$, where $z_{R}$ is a singularity of $X_{A_{t_R}}^s\cap \mathbb{C}^{*I}\cap B_R$. Then $(0,0)$ is in the closure of the set
$$W =\{(t,z) \in D\times \mathbb{C}^{*I} : z \in X_{A_t}^s \textrm{ and } z \textrm{ is a singular point of } X_{A_t}^s\}.$$

Then, by the Curve Selection Lemma \cite{M1}, there exists an analytic curve
$$\begin{array}{cccc}
p: & [0,\epsilon ] & \to & W \\
   & q & \mapsto & (t(q),z(q))
\end{array}$$
such that $p(q) = (t(q),z_1(q),\dots ,z_r(q),0,\dots,0)$, for all $q\neq 0$, and $p(0)=(0,0)$. For $1\leq i\leq r$, consider the Taylor expansions
\[
t(q) = t_0q^{\omega}+\cdots, \qquad z_i(q) = a_i q^{w_i} +\cdots,
\]
where $t_0,a_i\neq 0$ and $\nu,w_i>0$, for $i=1,\dots ,r$. Here, the three centered dots stand for
the higher order terms. Choose $a=(a_1,\dots,a_r,0,\dots,0) \in
\mathbb{C}^{*I}$ and $w=(w_1,\dots,w_r,0,\dots ,0) \in \mathbb{N}^m\setminus
\{0\}$. Consider the face $\Gamma_j$ of $(\Delta_j^{t(q)})^I=(\Delta_j^0)^I$
defined as the set where the map
\[
\begin{array}{cccc}
l_w^j : & (\Delta_j^0)^I & \to & \mathbb{R}_{+} \\
        & x:=(x_1,\dots,x_r,0,\dots ,0) & \mapsto & \displaystyle \sum_{i=1}^r x_iw_i
\end{array}
\]
takes its minimal value $d_j$ and such that $\sum_{j=1}^k\Gamma_j$ is a bounded face of $\sum_{j=1}^k\Delta_j$.

\begin{claim}\label{claim1}
	The point $a$ belongs to $X_{A_0^{*I}|_{\Gamma}}^s =
	((a_{ij}^{*I})_{0}|_{\Gamma_j})$.
\end{claim}

\begin{claimproof}
	Firstly, consider the $n\times k$ matrix
\[
(A_{t(q)}^{*I})|_{\Gamma}(z(q))  =
((a_{ij}^{*I})_{t(q)}|_{\Gamma_j}(z(q)).
\]

Now, note that

\begin{equation}\label{eachentry}
(a_{ij})_t(z) = ((a_{ij})_t)|_{\Gamma_j}(z) + \displaystyle \sum_{\alpha \notin \Gamma_j} c_{\alpha}z^{\alpha}.
\end{equation}

Consider a monomial component
\begin{equation}\label{monomial}
(c_{\alpha}z^{\alpha})|_{\Gamma_j} = c_{\alpha}z_1^{\alpha_1}\dots
z_r^{\alpha_r}
\end{equation}
of the face function $a_{ij}|_{\Gamma_j}$. Then over the curve $p$ we have
\[
\begin{aligned}
(c_{\alpha}z(q)^{\alpha})|_{\Gamma_j} & =
c_{\alpha}(a_1q^{w_1}+\cdots)^{\alpha_1}\ldots (a_rq^{w_r}+\cdots)^{\alpha_r}\\
 & =   c_{\alpha}a_1^{\alpha_1}\ldots a_r^{\alpha_r} q^{d_j}+\cdots.
\end{aligned}
\]
Therefore,
\begin{equation}\label{facefunction}
(a_{ij}^{*I})_{t(q)}|_{\Gamma_j}(z(q)) =
(a_{ij}^{*I})_{t(q)}|_{\Gamma_j}(a)q^{d_j}+\cdots
\end{equation}
and \[
(A_{t(q)}^{*I})|_{\Gamma}(z(q))  =
 ((a_{ij}^{*I})_{t(q)}|_{\Gamma_j}(a)q^{d_j}+\cdots ).
\]
Consider the set
	$$\mathcal{C} = \{(\mathcal{I},\mathcal{J}): \mathcal{I}\subset \{1,\dots ,n\} \textrm{, } \mathcal{J}\subset \{1,\dots ,k\} \textrm{ and } |\mathcal{I}|=|\mathcal{J}|=s \}.
	$$
Since the Newton polyhedron of $a_{ij}^{*I}$ is $\Delta_j^I$ for all $i=1,\dots ,n$ and $z(q)$ belongs to the determinantal variety $X_{A_{t(q)}^{*I}|_{\Gamma}}^s$,
which is defined by the $s$ size minors of $A_{t(q)}^{*I}|_{\Gamma}$, for each $(\mathcal{I},\mathcal{J})\in \mathcal{C}$, we have the zero polynomial
$$\begin{aligned}\det (((a_{ij}^{*I})_{t(q)}|_{\Gamma_j}(a)&q^{d_j}+\cdots )_{i\in\mathcal{I},\;j\in \mathcal{J}})\\ & = \det(((a_{ij}^{*I})_{t(q)}|_{\Gamma_j}(a))_{i\in\mathcal{I},\;j\in \mathcal{J}}) q^{\sum_{j \in \mathcal{J}}d_j}+\cdots = 0,\end{aligned}$$
which implies that all the coefficients of this polynomial are equal to zero, in particular
\[\det
(((a_{ij}^{*I})_{t(q)}|_{\Gamma_j}(a))_{i\in\mathcal{I},\;j\in
\mathcal{J}})=0.
\]
Taking $q\rightarrow 0$, for each $(\mathcal{I},\mathcal{J})\in \mathcal{C}$,
then $t(q)\rightarrow 0$. Therefore, the point $a$ belongs to the determinantal
variety $ X_{A_0^{*I}|_{\Gamma}}^s$.
\end{claimproof}

\begin{claim}\label{claim2}
The point $a$ is a singularity of $X_{A_0^{*I}|_{\Gamma}}^s$.
\end{claim}

\begin{claimproof}
We start by taking the partial derivative of \Eqref{monomial}, which gives us the following equation
\[
\frac{\partial }{\partial z_l}(c_{\alpha}z^{\alpha}) = \alpha_l c_{\alpha}
z_1^{\alpha_1}\dots z_l^{\alpha_l -1}\dots z_r^{\alpha_r}.
\]
Then, over the curve $p$
\[
\begin{aligned}
\frac{\partial}{\partial
z_l}(c_{\alpha}z(q)^{\alpha})|_{\Gamma_j} & =
\alpha_lc_{\alpha}(a_1q^{w_1}+\cdots )^{\alpha_1}\cdots (a_lq^{w_l}+\cdots
)^{\alpha_l -1}\cdots (a_rq^{w_r}+\cdots)^{\alpha_r}\\
  & =  \alpha_la_1^{\alpha_1}\cdots a_l^{\alpha_l-1}\cdots a_r^{\alpha_r}
  q^{d_j-w_l}+\cdots \\
  & =  \frac{\partial}{\partial
  z_l}(c_{\alpha}z(q)^{\alpha})|_{\Gamma_j}(a)\cdot q^{d_j-w_l}+\cdots.
\end{aligned}
\]

Thus, the following equation holds
\begin{equation} \label{derivative}
\frac{\partial}{\partial
z_l}((a_{ij}^{*I})_{t(q)}|_{\Gamma_j}) (z(q)) =
\displaystyle\frac{\partial}{\partial
z_l}((a_{ij}^{*I})_{t(q)}|_{\Gamma_j}) (a)\cdot
q^{d_j-w_l}+\cdots.
\end{equation}

It follows from \Eqref{eachentry} that
\[
((a_{ij}^{*I})_{t(q)}) (z(q)) =
((a_{ij}^{*I})_{t(q)}|_{\Gamma_j}) (z(q)) + \displaystyle
\sum_{\alpha \notin \Gamma_j} c_{\alpha}z(q)^{\alpha}
\]
and
\[
\displaystyle\frac{\partial}{\partial
z_l}((a_{ij}^{*I})_{t(q)}) (z(q)) =
\displaystyle\frac{\partial}{\partial
z_l}((a_{ij}^{*I})_{t(q)}|_{\Gamma_j}) (z(q)) +
\displaystyle\frac{\partial}{\partial z_l} (\displaystyle \sum_{\alpha
\notin \Gamma_j} c_{\alpha}z(q)^{\alpha}).
\]

Furthermore, by \Eqref{facefunction}, \Eqref{derivative} and the derivative product law, we have
	\begin{equation*}
		\begin{aligned}
			\frac{\partial}{\partial z_l} (\det ((a_{ij}^{*I})_{t(q)}|_{\Gamma_j}))_{i\in\mathcal{I},\;j\in
				\mathcal{J}} &(z(q)) \\
		 &=\frac{\partial}{\partial z_l} (\det
		 ((a_{ij}^{*I})_{t(q)}|_{\Gamma_j}))_{i\in\mathcal{I},\;j\in
		 	\mathcal{J}}(a) \cdot q^{(\sum_{j\in \mathcal{J}}d_j)-w_l}+ \cdots.
		\end{aligned}
		\end{equation*}

Then the jacobian matrix of the map given by the $s$ size minors of the
matrix $((a_{ij}^{*I})_{t(q)}|_{\Gamma_j}(a)q^{d_j}+\cdots
)_{i\in\mathcal{I},\;j\in \mathcal{J}}$ is the $\binom{k}{s}\binom{n}{s}
\times m$ matrix
\[
(\frac{\partial}{\partial z_l}(\det
((a_{ij}^{*I})_{t(q)}|_{\Gamma_j})_{i\in\mathcal{I},\;j\in
\mathcal{J}})(a)\cdot q^{(\sum_{j\in
\mathcal{J}}d_j)-w_l}+\cdots)_{(\mathcal{I},
\mathcal{J})\in \mathcal{C}}.
\]
Since $z(q)$ is a singularity of the determinantal variety
$X_{A_{t(q)}^{*I}|_{\Gamma}}^s$, this matrix has rank less than $(n-s+1)(k-s+1)$. This means that the $(n-s+1)(k-s+1)$ size minors of this jacobian matrix are all equal to zero.

Note that, for each $\theta=1,\dots ,(n-s+1)(k-s+1)$, the zero polynomial

\[\begin{aligned}
&(\det(\frac{\partial}{\partial z_{l_{\theta}}}(\det
((a_{ij}^{*I})_{t(q)}|_{\Gamma_{j}})_{i\in\mathcal{I}_{\theta},\;j\in
\mathcal{J}_{\theta}})(a) \cdot q^{(\sum_{j\in
\mathcal{J}_{\theta}}d_j)-w_{l_{\theta}}}+\cdots)_{(\mathcal{I}_{\theta},
\mathcal{J}_{\theta})\in \mathcal{C}})_{\theta} \\ &=
(\det(\frac{\partial}{\partial z_{l_{\theta}}}(\det
((a_{ij}^{*I})_{t(q)}|_{\Gamma_j})_{i\in\mathcal{I}_{\theta},\;j\in
\mathcal{J}_{\theta}})(a)_{\mathcal{I}_{\theta},
\mathcal{J}_{\theta}\in \mathcal{C}})_{\theta})\\& \times q^{(\sum_{\theta=1}^{(n-s+1)(k-s+1)}\sum_{j\in
\mathcal{J}_{\theta}}d_j)- \sum_{l_{\theta}=1}^{(n-s+1)(k-s+1)}w_{l_{\theta}}}+\cdots = 0,
\end{aligned}
\]
implies that all of its coefficients are zero, specifically
\[
\det(\frac{\partial}{\partial z_{l_{\theta}}}(\det
((a_{ij}^{*I})_{t(q)}|_{\Gamma_j})_{i\in\mathcal{I}_{\theta},\;j\in
\mathcal{J}_{\theta}})(a)_{\mathcal{I}_{\theta},
\mathcal{J}_{\theta}\in \mathcal{C}})_{\theta} = 0.
\]
Therefore, the jacobian matrix whose columns are
\[
\frac{\partial}{\partial
z_l}\det
((a_{ij}^{*I})_{t(q)}|_{\Gamma_j})_{i\in\mathcal{I},\;j\in
\mathcal{J}}(a)
\]
has also rank less than $(n-s+1)(k-s+1)$, where $l=1, \dots, m$.

Taking again $q\rightarrow 0$, then $t(q)\rightarrow 0$. Therefore, the point $a$ is a singularity of the determinantal variety $ X_{A_0^{*I}|_{\Gamma}}^s$.
\end{claimproof}

By Claims \ref{claim1} and \ref{claim2}, the matrix $ A_0^{*I}$ is not Newton
non-degenerate which contradicts Lemma \ref{lemmaadmissible}.

To prove the transversality, we also argue by contradiction. Suppose that there exists a sequence $\{(t_R,z_R)\}$ of points in $X_{\mathcal{A}}^s\cap (D\times \mathbb{C}^{*I})$ converging to $(0,0)$ and such that $X_{A_{t_R}}^s\cap \mathbb{C}^{*I}$ does not intersect the sphere $S_{||z_R||}$ transversally at $z_R$. Thus
	$$(T_{z_R}S_{||z_{R}||})^{\bot}  \subseteq (T_{z_R}(X_{A_{t_R}}^s\cap 		\mathbb{C}^{*I}))^{\bot}.
	$$
The orthogonal space $(T_{z_R}(X_{A_{t_R}}^s\cap \mathbb{C}^{*I}))^{\bot}$ is generated by the set
	$$G = \{grad (\det ((a_{ij}^{*I})_{t_R}(z_R))_{i\in \mathcal{I}, j\in \mathcal{J}}):(\mathcal{I}, \mathcal{J}) \in \mathcal{C}\}.
	$$
\noindent where the gradient of a function $f$ is $grad (f (z)) = \left(\overline{\frac{\partial f }{\partial z_1}(z)}, \dots, \overline{\frac{\partial f }{\partial z_m}(z)} \right)$ and $\overline{\frac{\partial f}{\partial z_i}(z)}$ denotes the complex conjugation of $\frac{\partial f}{\partial z_i}(z)$, $i=1, \dots , m$. For simplification purposes, we will denote $\overline{\frac{\partial f}{\partial z_i}(z)} = \frac{\overline{\partial} f}{\partial z_i}(z)$, $i=1, \dots , m$.

In general, the set $G$ is not linearly independent, then $G$ is not a basis for $(T_{z_R}(X_{A_{t_R}}^s\cap \mathbb{C}^{*I}))^{\bot}$. Our first task is finding, if necessary, a subsequence $(t_{R_{\gamma}},z_{R_{\gamma}})$ of $(t_R,z_R)$, such that
$$G_{\gamma}=\{grad (\det ((a_{ij}^{*I})_{t_{R_{\gamma}}}(z_{R_{\gamma}}))_{i\in \mathcal{I}, j\in \mathcal{J}}):(\mathcal{I}, \mathcal{J}) \in \mathcal{C}_{\gamma}\}$$
forms a basis for $(T_{z_{R_{\gamma}}}(X_{A_{t_{R_{\gamma}}}}^s\cap \mathbb{C}^{*I}))^{\bot}$.

Firstly, we take subsets $\mathcal{P}_{\gamma}$ of $\mathbb{N}$ such that
$$G_{\gamma} = \{grad (\det ((a_{ij}^{*I})_{t_{R_{\gamma}}}(z_{R_{\gamma}}))_{i\in \mathcal{I}, j\in \mathcal{J}}):(\mathcal{I}, \mathcal{J}) \in \mathcal{C}_{\gamma}\}$$
is a basis for $(T_{z_{R_{\gamma}}}(X_{A_{t_{R_{\gamma}}}}^s\cap \mathbb{C}^{*I}))^{\bot}$, where $\mathcal{C}_{\gamma}\subset \mathcal{C}$ and $1/R_{\gamma} \in \mathcal{P}_{\gamma}$.

Since there are only finitely many subsets $\mathcal{P}_{\gamma}$, there exists $\gamma_{0}$ such that $\mathcal{P}_{\gamma_{0}}$ is not finite. Therefore, $\{(t_{R_{\gamma_{0}}},z_{R_{\gamma_{0}}})\}_{\frac{1}{R_{\gamma_{0}}} \in \mathcal{P}_{\gamma_{0}}}$ is the desired subsequence. To simplify, we denote this subsequence by $\{(t_{R_{\gamma}},z_{R_{\gamma}})\}$.  Then we can write the vector $z_{R_{\gamma}}\in (T_{z_{R_{\gamma}}}S_{||z_{R_{\gamma}}||})^{\bot}$ as a linear combination of vectors from $G_{\gamma}$, \ie there exist $\lambda_{\mathcal{I},\mathcal{J}}$ satisfying
$$z_{R_{\gamma}}= \sum_{(\mathcal{I}, \mathcal{J})\in \mathcal{C}_{\gamma}} \lambda_{\mathcal{I},\mathcal{J}} \cdot grad(\det ((a_{ij}^{*I})_{t_{R_{\gamma}}}(z_{R_{\gamma}}))_{i\in \mathcal{I}, j\in \mathcal{J}}).
$$

We observe that some of the coefficients $\lambda_{\mathcal{I},\mathcal{J}}$ in the above linear combination may be zero. However, by the same arguments, if necessary, we can consider a subsequence $(t_{R_{\gamma_{\alpha}}},z_{R_{\gamma_{\alpha}}} )$ of $(t_{R_{\gamma}},z_{R_{\gamma}} )$ given by
$$
z_{R_{\gamma_{\alpha}}} = \sum_{(\mathcal{I},\mathcal{J})	\in \mathcal{C}_{\gamma_{\alpha}}} \lambda_{\mathcal{I},\mathcal{J}} \cdot grad(\det ((a_{ij}^{*I})_{t_{R_{\gamma_{\alpha}}}} (z_{R_{\gamma_{\alpha}}}))_{i\in \mathcal{I}, j\in \mathcal{J}}),
$$
where $\lambda_{\mathcal{I}, \mathcal{J}} \neq 0$, for all $(\mathcal{I}, \mathcal{J}) \in \mathcal{C}_{\gamma_{\alpha}}$ and $1/R_{\gamma_{\alpha}} \in \mathcal{P}_{\gamma_{\alpha}}$. We will denote this subsequence simply by $\{(t_{R_{\alpha}},z_{R_{\alpha}})\}$ and we denote by $\mathcal{C}_{\alpha}$ the set of $(\mathcal{I},\mathcal{J})\in \mathcal{C}$ such that
$$\{grad(\det ((a_{ij}^{*I})_{t_{R_{\alpha}}}(z_{R_{\alpha}}))_{i\in \mathcal{I}, j\in \mathcal{J}})\}$$ is linearly independent and $\lambda_{\mathcal{I},\mathcal{J}} \neq 0$ and $1/\alpha \in \mathcal{P}_{\alpha}$.

Since the subsequence  $(t_{R_{\alpha}},z_{R_{\alpha}})$ also converges to $(0,0)$, the point $(0,0)$ belongs to the closure of the set consisting of points $(t,z)\in D\times \mathbb{C}^{*I}$ such that
	$$z\in X_{A_t}^s\; \textrm{and}\; 	z = \sum_{(\mathcal{I},\mathcal{J})	\in \mathcal{C}_{\alpha}} \lambda_{\mathcal{I},\mathcal{J}} \cdot grad(\det ((a_{ij}^{*I})_t(z))_{i\in \mathcal{I}, j\in \mathcal{J}}).
	$$
By the Curve Selection Lemma \cite{M1}, there exists a real analytic curve
$$(t(q),z(q)) = (t(q), z_1(q),\dots ,z_r(q),0,\dots ,0)$$
and Laurent series $\lambda_{\mathcal{I},\mathcal{J}}$, $(\mathcal{I},\mathcal{J})\in \mathcal{C}_{\alpha}$, such that
\begin{enumerate}
\item $(t(0),z(0)) = (0,0)$;
\item $(t(q),z(q)) \in D\times \mathbb{C}^{*I}$, for $q\neq 0$;
\item $z(q)\in X_{A_{t(q)}}^s$;
\item $z(q) = \sum_{(\mathcal{I},\mathcal{J})	\in \mathcal{C}_{\alpha}} \lambda_{\mathcal{I},\mathcal{J}} \cdot grad(\det (((a_{ij}^{*I})_{t(q)}(z(q)))_{i\in \mathcal{I}, j\in \mathcal{J}})).$
\end{enumerate}

Consider the Taylor expansions
\[
t(q) = t_0q^{\omega}+\cdots, \qquad z_i(q) = a_i q^{w_i} +\cdots,
\]
where $t_0,a_i\neq 0$ and $\nu,w_i>0$, for $i=1, \dots, r$. Consider also the Laurent expansions
	$$\lambda_{\mathcal{I},\mathcal{J}}(q) = \beta_{\mathcal{I},\mathcal{J}}\cdot q^{u_{\mathcal{I},\mathcal{J}}}+\cdots,
	$$
where $\beta_{\mathcal{I},\mathcal{J}}\neq 0$.

We define $a=(a_1,\dots ,a_r, 0,\dots ,0)\in \mathbb{C}^{*I}$, $w = (w_1,\dots ,w_r,0,\dots ,0)$ and $d_j$ the minimum value of the function $l_w^j:(\Delta_j^{t(q)})^I\to \mathbb{R}_{+}$, with $\Gamma_j$ being the face of $(\Delta_j^{t(q)})^I=(\Delta_j^0)^I$ where this function takes this minimum value and such that $\sum_{j=1}^k\Gamma_j$ is a bounded face of $\sum_{j=1}^k \Delta_k$.

\begin{claim}\label{claim3}
	There exists $\tilde{\mathcal{C}} \subset \mathcal{C}_{\alpha}$ such that
\[
\sum_{(\mathcal{I},\mathcal{J})\in \tilde{\mathcal{C}}} \beta_{\mathcal{I},\mathcal{J}} \sum_{l=1}^r w_l \overline{a_l} \frac{\overline{\partial}}{\partial z_l}\det ((a_{ij}^{*I})_0|_{\Gamma_j})_{i\in \mathcal{I}, j\in \mathcal{J}}(a)  \neq 0 .
\]
\end{claim}

\begin{claimproof}
Indeed, since
\[
\begin{aligned}
& grad(\det ((a_{ij}^{*I})_{t(q)}|_{\Gamma_j}(z(q)))_{i\in \mathcal{I}, j\in \mathcal{J}}) =\\ &(\frac{\overline{\partial}}{\partial z_1}\det ((a_{ij}^{*I})_{t(q)}|_{\Gamma_j})_{i\in \mathcal{I}, j\in \mathcal{J}}(a)\cdot q^{(\sum_{j\in \mathcal{J}}d_j) -w_1}+\cdots, \cdots ,\\ &\frac{\overline{\partial}}{\partial z_r}\det ((a_{ij}^{*I})_{t(q)}|_{\Gamma_j})_{i\in \mathcal{I}, j\in \mathcal{J}}(a)\cdot q^{(\sum_{j\in \mathcal{J}}d_j) -w_r}+\cdots,0,\cdots, 0 ),
\end{aligned}
\]
by $(iv)$, we have
\[
\begin{aligned}
& a_lq^{w_l}+\cdots = z_l(q) = \\ &\sum_{(\mathcal{I},\mathcal{J})\in \mathcal{C}_{\alpha}} \beta_{\mathcal{I},\mathcal{J}} \frac{\overline{\partial}}{\partial z_l}\det ((a_{ij}^{*I})_{t(q)}|_{\Gamma_j})_{i\in \mathcal{I}, j\in \mathcal{J}}(a)\cdot q^{\sum_{j\in \mathcal{J}}d_j +u_{\mathcal{I},\mathcal{J}} -w_l}+\cdots,
\end{aligned}
\]
for $l=1,\dots ,m$.

We choose the set $\tilde{\mathcal{C}} \subset \mathcal{C}_{\alpha}$, such that, for each $(\mathcal{I},\mathcal{J}), (\tilde{\mathcal{I}},\tilde{\mathcal{J}})$ in $\tilde{\mathcal{C}}$, we have
\[
\begin{aligned}
(\sum_{j\in \mathcal{J}} d_j) + u_{\mathcal{I},\mathcal{J}} = (\sum_{j\in \tilde{\mathcal{J}}} d_j) &+ u_{\tilde{\mathcal{I}},\tilde{\mathcal{J}}} =\\ & \min \{(\sum_{j\in \mathcal{J}} d_j) + u_{\mathcal{I},\mathcal{J}} : (\mathcal{I},\mathcal{J}) \in \mathcal{C}_{\alpha}\}.
\end{aligned}
\]
Then, $w_l = (\sum_{j \in \mathcal{J}} d_j)+u_{\mathcal{I},\mathcal{J}} -w_l$, for all $(\mathcal{I},\mathcal{J}) \in \tilde{\mathcal{C}}$.

We may reorder $w_1,\dots ,w_r$, if necessary, such that $w_1 =\cdots =w_b< w_c$ $(b<c\leq r)$. Therefore,
$$\sum_{(\mathcal{I},\mathcal{J}) \in \tilde{\mathcal{C}}} \beta_{\mathcal{I},\mathcal{J}} \frac{\overline{\partial}}{\partial z_l}\det ((a_{ij}^{*I})_{t(q)}|_{\Gamma_j})_{i\in \mathcal{I}, j\in \mathcal{J}}(a) = \left\{ \begin{array}{cl} a_l, & 1\leq l\leq b, \\ 0, & b<l\leq r \end{array} \right. . $$

Multiplying both sides of the last equality by $w_l\overline{a_l}$, we have
$$
\sum_{(\mathcal{I},\mathcal{J}) \in \tilde{\mathcal{C}}} \beta_{\mathcal{I},\mathcal{J}} w_l \overline{a_l}\frac{ \overline{\partial}}{\partial z_l}\det ((a_{ij}^{*I})_{t(q)}|_{\Gamma_j})_{i\in \mathcal{I}, j\in \mathcal{J}}(a) = \left\{ \begin{array}{cl} w_l |a_l|^2, & 1\leq l\leq b, \\ 0, & b<l\leq r \end{array} \right. .
$$

Taking a sum over $1\leq l\leq r$ and taking $q\to 0$ we have
\[
\sum_{(\mathcal{I},\mathcal{J}) \in \tilde{\mathcal{C}}} \beta_{\mathcal{I},\mathcal{J}} \sum_{l=1}^r w_l \overline{a_l} \frac{\overline{\partial}}{\partial z_l}\det ((a_{ij}^{*I})_0|_{\Gamma_j})_{i\in \mathcal{I}, j\in \mathcal{J}}(a) = \sum_{l=1}^r w_l |a_l|^2 \neq 0 .
\]
\end{claimproof}

\begin{claim}\label{claim4}
	For each $(\mathcal{I},\mathcal{J}) \in \tilde{\mathcal{C}}$, the following equality holds $$ \sum_{l=1}^r w_l a_l
	\frac{\partial}{\partial z_l}\det ((a_{ij}^{*I})_0|_{\Gamma_j})_{i\in
	\mathcal{I}, j\in \mathcal{J}}(a) = 0.$$
\end{claim}

\begin{claimproof}
The polynomial $\det ((a_{ij}^{*I})_{t(q)}|_{\Gamma_j}(a))_{i\in \mathcal{I}, j\in \mathcal{J}}$ is weighted homogeneous with respect to the weight $w$ and it has weighted degree $$\sum_{j\in \mathcal{J}} d_j,$$ it follows from the Euler identity that
$$(\sum_{j\in \mathcal{J}} d_j)\cdot \det ((a_{ij}^{*I})_{t(q)}|_{\Gamma_j}(a))_{i\in \mathcal{I}, j\in \mathcal{J}} = \sum_{l=1}^r w_l a_l \frac{\partial}{\partial z_l}\det ((a_{ij}^{*I})_{t(q)}|_{\Gamma_j})_{i\in \mathcal{I}, j\in \mathcal{J}}(a).
$$

Taking $q\to 0$, by the same arguments of Claim 1, $a\in X_{A_0^{*I}}^s$. Then
\[
\det ((a_{ij}^{*I})_0|_{\Gamma_j}(a))_{i\in \mathcal{I}, j\in \mathcal{J}}=0.
\]
\end{claimproof}

Combining Claim \ref{claim3} and Claim \ref{claim4}, we get the contradiction

$$0 = \sum_{(\mathcal{I},\mathcal{J}) \in \tilde{\mathcal{C}}} \beta_{\mathcal{I},\mathcal{J}}
\sum_{l=1}^r w_l \overline{a_l} \frac{\overline{\partial}}{\partial z_l}\det
((a_{ij}^{*I})_0|_{\Gamma_j})_{i\in \mathcal{I}, j\in
\mathcal{J}}(a) \neq 0,$$ 
and the result follows.
\end{proof}

\begin{corollary} \label{cororadius}
In addition to the conditions of Proposition \ref{radius}, if the Newton polyhedra $\Delta_j^t$ of $(a_{ij})_t$ are convenient, then there exists a positive number $R > 0$ such that for any $t$ sufficiently small, the set $X_{A_t}^s \cap B_R$ is smooth outside the origin and intersects transversely with $S_r$ for any $r<R$, where $B_R$ (respectively, $S_r$) is the open ball (respectively, the sphere) with center at the origin $0 \in \mathbb{C}^m$ and radius $R$ (respectively, $r$).
\end{corollary}

\begin{corollary}\label{vanishing-constant}
If for all t sufficiently small, the following two conditions are satisfied:
\begin{enumerate}
\item the matrix $A_t=((a_{ij})_t)$, is Newton non-degenerate;
\item the Newton polyhedron $\Delta_j^t$ of $(a_{ij})_t$ is convenient and independent of $t$, for all $j=1,\dots ,k$.
\end{enumerate}
Then the vanishing Euler characteristic of $X_{A_t}^s$ is independent of $t$, \ie
$$\nu (X_{A_t}^s,0) = \nu (X_{A_0}^s,0).$$
\end{corollary}
\begin{proof}
It follows directly from Corollaries \ref{corovanishing} and \ref{cororadius}.
\end{proof}

Now, given a determinantal deformation of $X_{A}^s$, consider the functions $f_{k+1},\dots ,f_p :(\mathbb{C}^m,0)\to (\mathbb{C},0)$ and function germs $F_{k+1},\dots ,F_p :(\mathbb{C}^m\times \mathbb{C},0)\to (\mathbb{C},0)$ such that $F_i(x,0)=f_i(x)$ for all $x\in \mathbb{C}^m$ and $i=k+1,\dots ,p$. For each $i=k+1,\dots ,p$, we fix small enough representatives $f_i :B_{\varepsilon} \to (\mathbb{C},0)$, where $B_{\varepsilon}$ is the open ball centered at the origin with radius $\varepsilon >0$, and we set $f_{(i,t)}(x) := F_i(x,t)$, $i=1,\dots ,p$. Therefore, we can consider the map germ $\tilde{\mathcal{A}}:=(\mathcal{A},F_{k+1},\dots ,F_p ): (\mathbb{C}^m\times \mathbb{C},0)\to (M_{n,k}\times\mathbb{C}^p,0)$ and we have a family of fibers
	$$\{X_{A_t}^s\cap f_{({k+1},t)}^{-1}(0)\cap \cdots \cap  f_{(p,t)}^{-1}(0)\}_{t\in D}.
	$$

The properties of fibers of functions defined on determinantal singularities were also studied, for instance, by Ament, Nu\~no-Ballesteros, Or\'efice-Okamoto and Tomazzela \cite{ABOT}, Carvalho, Nu\~no-Ballesteros, \linebreak Or\'efice-Okamoto and Tomazzela \cite{CBBT} and Menegon and Pereira \cite{AMP}.

From now on, we will denote by $\Delta_1^t, \dots, \Delta_k^t, \Delta_{k+1}^t, \dots, \Delta_{p}^t$ the Newton polyhedra of the columns of the matrix $A_t=((a_{ij})_t)$ and the Newton polyhedra of the functions $f_{(k+1,t)},\dots ,f_{(p,t)}$, respectively. With this notation, we introduce the following result.

\begin{proposition}
\label{radius-section}
Suppose that for all t sufficiently small, the following two conditions are satisfied:
\begin{enumerate}
\item the variety $X_{A_t}^s \cap f_{(k+1,t)}^{-1}(0)\cap \cdots \cap f_{(p,t)}^{-1}(0)$ is Newton non-degenerate, where $A_t=((a_{ij})_t)$;
\item the Newton polyhedra $\Delta_1^t, \dots, \Delta_k^{t}, \Delta_{k+1}^t,\dots ,\Delta_p^t$ are independent of $t$.
\end{enumerate}
Then there exists a positive number $R > 0$ such that for any admissible coordinate subspace $\mathbb{C}^I$ of $A_0$, $f_{(k+1,0)}, \dots ,f_{(p,0)}$ and any $t$ sufficiently small  the variety $$X_{A_t}^s \cap f_{(k+1,t)}^{-1}(0)\cap \cdots \cap f_{(p,t)}^{-1}(0) \cap \mathbb{C}^{*I} \cap B_R$$ is non-singular and intersects transversely with $S_r$ for any $r<R$, where $B_R$ (respectively, $S_r$) is the open ball (respectively, the sphere) with center at the origin $0 \in \mathbb{C}^m$ and radius $R$ (respectively, $r$).
\end{proposition}
\begin{proof}
The proof follows from the proof of Proposition \ref{radius} and \cite[Proposition 3.1]{EyralOka}.
\end{proof}

\begin{corollary}
\label{cororadius-section}
In addition to the conditions of Proposition \ref{radius-section}, if the Newton polyhedra $\Delta_1^t, \dots, \Delta_k^t, \Delta_{k+1}^t,\dots ,\Delta_p^t$ are convenient, then there exists a positive number $R > 0$ such that for any $t$ sufficiently small, the variety $X_{A_t}^s \cap f_{k+1}^{-1}(0)\cap \cdots \cap f_p^{-1}(0) \cap B_R$ is smooth outside the origin and intersects transversely with $S_r$ for any $r<R$, where $B_R$ (respectively, $S_r$) is the open ball (respectively, the sphere) with center at the origin $0 \in \mathbb{C}^m$ and radius $R$ (respectively, $r$).
\end{corollary}


Let $h:(\mathbb{C}^m,0)\to (\mathbb{C},0)$ be a generic linear form with respect to $X_A^s$. Then, the variety $X_A^s\cap h^{-1}(0)$ is Newton non-degenerate and we can suppose that the Newton polyhedron $\Delta_h$ of $h$ is convenient. Moreover, if $X_{A_t}^s$ is a determinantal deformation of an IDS $X_A^s$,  we can consider a deformation $X_{A_t}^s \cap h_t^{-1}(0)$ of $X_A^s \cap h^{-1}(0)$ where for each $t$, $h_t: (\mathbb{C}^m,0)\to (\mathbb{C},0) $ is a generic linear form with respect to $X_{A_t}^s$ such that $\Delta_{h_t}$ is convenient for all $t$. Therefore, we have the following corollary.

\begin{corollary}\label{toppolar-constant}
Suppose that for all t sufficiently small, the following conditions are satisfied:
\begin{enumerate}
\item the variety $X_{A_t}^s$ is Newton non-degenerate;
\item the Newton polyhedra $\Delta_j^t$ of $(a_{ij})_t$ are convenient and independent of $t$, for all $j=1,\dots ,k$.
\end{enumerate}
Then the vanishing Euler characteristic of $X_{A_t}^s \cap h_t^{-1}(0)$ and consequently the top polar multiplicity of $X_{A_t}^s$ are independent of $t$, \ie
$$\nu (X_{A_t}^s \cap h_{t}^{-1}(0),0) = \nu (X_{A_0}^s \cap h_{0}^{-1}(0),0)$$
and
$$m_d (X_{A_t}^s,0) = m_d(X_{A_0}^s,0).$$
\end{corollary}

\begin{proof}
The first equality follows directly from Corollaries \ref{corovanishing} and \ref{cororadius-section}. The second equation follows directly from the first equality, \Eqref{toppolar2} and Corollary \ref{vanishing-constant}.
\end{proof}

Finally, we have all the tools necessary to prove Theorem \ref{equisingularity}.

\begin{proof}[Proof of Theorem \ref{equisingularity}]
Firstly, since the Newton polyhedra $\Delta_j^t$ are convenient and independent of $t$ for each $j=1,\dots,k$ and the matrix $A_t$ is Newton non-degenerate, then by Corollary \ref{cororadius}, there exists a positive number $R$ such that for any $t$ sufficiently small, the set $X_{A_t}^s\cap B_{R}$ is smooth outside the origin, where $B_{R}$ is the open ball with center at the origin and radius $R$. Therefore, this family is good.

Moreover, by Corollary \ref{toppolar-constant},
\begin{equation}
\label{toppolar0}
m_d(X_{A_t}^s,0)= m_d(X_{A_0}^s,0).
\end{equation}

Now, considering $h_1, \dots, h_l$ generic linear forms and applying successively the  Corollary \ref{toppolar-constant}
we have the following relation,
$$\nu(X_{A_t}^s\cap h_{(t,1)}^{-1}(0) \cap \cdots \cap h_{(t,l)}^{-1}(0),0)=\nu(X_{A_0}^s \cap h_{(0,1)}^{-1}(0) \cap \cdots \cap h_{(0,l)}^{-1}(0),0),$$
where $h_{(t,i)}$ are generic linear forms with respect to
	$$X_{A_t}^s\cap h_{(t,1)}^{-1}(0) \cap \cdots \cap h_{(t,i-1)}^{-1}(0)
	$$
for $i=1,\dots,l$. Hence,
\begin{equation*}
\label{toppolarl}
m_d(X_{A_t}^s\cap h_{(t,1)}^{-1}(0) \cap \cdots \cap h_{(t,l)}^{-1}(0),0)=m_d((X_{A_0}^s \cap h_{(0,1)}^{-1}(0) \cap \cdots \cap h_{(0,l)}^{-1}(0),0).
\end{equation*}

In addition, applying successively \cite[Lemma $2.6$]{GGR}, we have
\begin{equation*}\label{polarcortes}
m_{d-l} (X_{A_t}^s\cap h_{(t,1)}^{-1}(0) \cap \cdots \cap h_{(t,l)}^{-1}(0),0)=m_{d-l}(X_{A_t}^s,0),
\end{equation*}
for $l=1,\dots,d$.
Combining both equations above, we have
\begin{equation*}
\begin{aligned}
m_{d-l} (X_{A_t}^s,0) & =  m_{d-l}(X_{A_t}^s\cap h_{(t,1)}^{-1}(0) \cap \cdots \cap h_{(t,l)}^{-1}(0),0)
\\
& =  m_{d-l}((X_{A_0}^s \cap h_{(0,1)}^{-1}(0) \cap \cdots \cap h_{(0,l)}^{-1}(0),0) \\
& =  m_{d-l} (X_{A_0}^s,0),
\end{aligned}
\end{equation*}
for all $l=1,\dots , d$.

Therefore, $m_j(X_{A_t}^s,0)=m_j(X_{A_0}^s,0)$, for all $0\leq j\leq d$. Hence, the family $\left\{(X_{A_t}^s, 0) \right\}_{ t\in D }$ is Whitney equisingular, by \cite[Theorem 5.3]{BBT4}.
\end{proof}

\begin{example}
Let $\left\{(X_{A_t}^2, 0) \right\}_{ t\in D }$ be the family of 2-dimensional determinantal singularities defined by the germ $A_t: (\mathbb{C}^4,0) \rightarrow (M_{2,3},0)$, given by the matrix
$$\scalemath{0.91}{
\left[
\begin{array}{ccc}
	2x + 2y +ty^2 + z - 3w & 2x + 3y + 2ty^2  + 2z - 5w & 3x +2y + ty^2  + 2z - 3w\\
	 3x +3y+ ty^2 + 2z - 4w & 3x +4y + 2ty^2 + 4z - 7w & 5x +3y + ty^2 + 3z - 3w
	\end{array}
	\right].
}
$$

For all $t\in D$, the matrix $A_t$ is Newton non-degenerate, $\Delta^t_{j}$ is convenient and independent of $t$, for all $j=1,\dots ,k$.
Then, by Theorem \ref{equisingularity}, $\left\{(X_{A_t}^2, 0) \right\}_{ t\in D}$ is Whitney equisingular.
\end{example}
\begin{remark}
The matrix $A_t$ in the last example is $\mathcal{G}$-equivalent to

$$\tilde{A}_t = \left[\begin{array}{ccc}
x & y +ty^2 & z\\
y & z & w
\end{array} \right].$$
Therefore, from the previous example, any family $\{(X_{A_t}^s,0)\}_{t\in D}$, where $A_t$ is $\mathcal{G}$-equivalent to $\tilde{A}_t$ is also Whitney equisingular, even though the Newton polyhedra $\Delta_j$ determined by the matrix $\tilde{A}_t$ are not convenient.
\end{remark}

\begin{center}
	{ \bf Acknowledgments}
\end{center}

We would like to thank Anne Fr\"uhbis-Kr\"uger, Bruna
Or\'efice-Okamoto, Juan Jos\'e Nu\~no-Ballesteros and Maria Aparecida Soares
Ruas
for all helpful discussions that they have provided along the
development of this work.

\bibliographystyle{amsalpha-lmp}
\bibliography{references}

\providecommand{\bysame}{\leavevmode\hbox to3em{\hrulefill}\thinspace}
\providecommand{\MR}{\relax\ifhmode\unskip\space\fi MR }
\providecommand{\MRhref}[2]{%
  \href{http://www.ams.org/mathscinet-getitem?mr=#1}{#2}
}
\providecommand{\href}[2]{#2}
\begin{thebibliography}{\textsc{CNnBOOT22}}

\bibitem[\textsc{ANnOOT16}]{ABOT}
\textsc{D.~A.~H. Ament}, \textsc{J.~J. Nu\~{n}o{-}Ballesteros},
  \textsc{B.~Or\'{e}fice-Okamoto}, and \textsc{J.~N. Tomazella}, \emph{The
  {E}uler obstruction of a function on a determinantal variety and on a curve},
  Bull. Braz. Math. Soc. (N.S.) \textbf{47} (2016), no.~3, 955--970.
  \MR{3549078}

\bibitem[\textsc{BrHe98}]{Bruns2}
\textsc{W.~Bruns} and \textsc{J.~Herzog}, \emph{Cohen-{M}acaulay rings},
  Revised edition, Cambridge University Press, New York (1998).

\bibitem[\textsc{BrTa04}]{BruceTari}
\textsc{J.~W. Bruce} and \textsc{F.~Tari}, \emph{On families of square
  matrices}, Proc. London Math. Soc. (3) \textbf{89} (2004), no.~3, 738--762.
  \MR{2107013}

\bibitem[\textsc{Bru03}]{Bruce}
\textsc{J.~W. Bruce}, \emph{On families of symmetric matrices}, vol.~3, 2003,
  Dedicated to Vladimir I. Arnold on the occasion of his 65th birthday,
  pp.~335--360, 741. \MR{2025264}

\bibitem[\textsc{BrVe98}]{Bruns1}
\textsc{W.~Bruns} and \textsc{U.~Vetter}, \emph{Determinantal rings},
  Springer-Verlang, New York (1998).

\bibitem[\textsc{CNnBOOT22}]{CBBT}
\textsc{R.~S. Carvalho}, \textsc{J.~J. Nu\~{n}o Ballesteros},
  \textsc{B.~Or\'{e}fice-Okamoto}, and \textsc{J.~N. Tomazella},
  \emph{Equisingularity of families of functions on isolated determinantal
  singularities}, Bull. Braz. Math. Soc. (N.S.) \textbf{53} (2022), no.~1,
  1--20. \MR{4379420}

\bibitem[\textsc{EbGZ09}]{EGZ}
\textsc{W.~Ebeling} and \textsc{S.~M. Gusein-Zade}, \emph{On indices of
  $1$-forms on determinantal singularities}, Tr. Mat. Inst. Steklova
  \textbf{267} (2009), 119--131.

\bibitem[\textsc{Est07}]{Esterov}
\textsc{A.~Esterov}, \emph{Determinantal singularities and {N}ewton polyhedra},
  Tr. Mat. Inst. Steklova \textbf{259} (2007), no.~Anal. i Osob. Ch. 2, 20--38.
  \MR{2433675}

\bibitem[\textsc{EyOk17}]{EyralOka}
\textsc{C.~Eyral} and \textsc{M.~Oka}, \emph{Non-compact {N}ewton boundary and
  {W}hitney equisingularity for non-isolated singularities}, Adv. Math.
  \textbf{316} (2017), 94--113. \MR{3672903}

\bibitem[\textsc{FK99}]{AFG1}
\textsc{A.~Fr\"{u}hbis-Kr\"{u}ger}, \emph{Classification of simple space curve
  singularities}, Comm. Algebra \textbf{27} (1999), no.~8, 3993--4013.
  \MR{1700205}

\bibitem[\textsc{FKNe10}]{AFG2}
\textsc{A.~Fr\"{u}hbis-Kr\"{u}ger} and \textsc{A.~Neumer}, \emph{Simple
  {C}ohen-{M}acaulay codimension 2 singularities}, Comm. Algebra \textbf{38}
  (2010), no.~2, 454--495. \MR{2598893}

\bibitem[\textsc{Gaf93}]{Gaffney}
\textsc{T.~Gaffney}, \emph{Polar multiplicities and equisingularity of map
  germs}, Topology \textbf{32} (1993), no.~1, 185--223. \MR{1204414}

\bibitem[\textsc{Gaf96}]{Gaffney2}
\bysame, \emph{Multiplicities and equisingularity of {ICIS} germs}, Invent.
  Math. \textbf{123} (1996), no.~2, 209--220. \MR{1374196}

\bibitem[\textsc{GGJR19}]{GGR}
\textsc{T.~Gaffney}, \textsc{N.~G. Grulha~Jr.}, and \textsc{M.~A.~S. Ruas},
  \emph{The local {E}uler obstruction and topology of the stabilization of
  associated determinantal varieties}, Math. Z. \textbf{291} (2019), no.~3-4,
  905--930. \MR{3936093}

\bibitem[\textsc{GrPf08}]{GP}
\textsc{G.-M. Greuel} and \textsc{G.~Pfister}, \emph{A {\bf {s}ingular}
  introduction to commutative algebra}, extended ed., Springer, Berlin, 2008,
  With contributions by Olaf Bachmann, Christoph Lossen and Hans
  Sch\"{o}nemann, With 1 CD-ROM (Windows, Macintosh and UNIX). \MR{2363237}

\bibitem[\textsc{Har77}]{Hartshorne}
\textsc{R.~Hartshorne}, \emph{Algebraic geometry}, Graduate Texts in
  Mathematics, No. 52, Springer-Verlag, New York-Heidelberg, 1977. \MR{463157}

\bibitem[\textsc{L{\^e}Te81}]{LT}
\textsc{D.~T. L{\^e}} and \textsc{B.~Teissier}, \emph{Vari\'et\'es polaires
  locales et classes de {C}hern des vari\'et\'es singuli\`eres}, Ann. of Math.
  (2) \textbf{114} (1981), no.~3, 457--491. \MR{634426}

\bibitem[\textsc{Mat12}]{Mather}
\textsc{J.~Mather}, \emph{Notes on topological stability}, Bull. Amer. Math.
  Soc. (N.S.) \textbf{49} (2012), no.~4, 475--506. \MR{2958928}

\bibitem[\textsc{MePe22}]{AMP}
\textsc{A.~Menegon} and \textsc{M.~S. Pereira}, \emph{On deformations of
  isolated singularity function}, arXiv:2206.10035v1 (2022).

\bibitem[\textsc{Mil68}]{M1}
\textsc{J.~Milnor}, \emph{Singular points of complex hypersurfaces}, Annals of
  Mathematics Studies, No. 61, Princeton University Press, Princeton, N.J.;
  University of Tokyo Press, Tokyo, 1968. \MR{0239612}

\bibitem[\textsc{NBOOT13}]{BBT}
\textsc{J.~J. Nu{\~n}o-Ballesteros}, \textsc{B.~Or\'efice-Okamoto}, and
  \textsc{J.~N. Tomazella}, \emph{The vanishing {E}uler characteristic of an
  isolated determinantal singularity}, Israel J. Math. \textbf{197} (2013),
  no.~1, 475--495. \MR{3096625}

\bibitem[\textsc{NnOOT18}]{BBT4}
\textsc{J.~J. Nu\~{n}o{-}Ballesteros}, \textsc{B.~Or\'{e}fice-Okamoto}, and
  \textsc{J.~N. Tomazella}, \emph{Equisingularity of families of isolated
  determinantal singularities}, Math. Z. \textbf{289} (2018), no.~3-4,
  1409--1425. \MR{3830255}

\bibitem[\textsc{Oka97}]{Oka}
\textsc{M.~Oka}, \emph{Non-degenerate complete intersection singularity},
  Actualit\'es Math\'ematiques. [Current Mathematical Topics], Hermann, Paris,
  1997. \MR{1483897}

\bibitem[\textsc{Oka10}]{Oka3}
\bysame, \emph{Non-degenerate mixed functions}, Kodai Math. J. \textbf{33}
  (2010), no.~1, 1--62. \MR{2732230}

\bibitem[\textsc{Per10}]{MP}
\textsc{M.~S. Pereira}, \emph{Variedades determinantais e singularidades de
  matrizes}, Ph.D. thesis, 2010.

\bibitem[\textsc{PeRu14}]{MC}
\textsc{M.~S. Pereira} and \textsc{M.~A.~S. Ruas}, \emph{Codimension two
  determinantal varieties with isolated singularities}, Math. Scand.
  \textbf{115} (2014), no.~2, 161--172. \MR{3291723}

\bibitem[\textsc{Tho69}]{Thom}
\textsc{R.~Thom}, \emph{Ensembles et morphismes stratifi\'{e}s}, Bull. Amer.
  Math. Soc. \textbf{75} (1969), 240--284. \MR{239613}

\end{thebibliography}

\end{document}